\documentclass[10pt]{siamltex}

\usepackage{graphicx}
\usepackage{epsfig}
\usepackage{amsmath}
\usepackage{amsfonts}
\usepackage{amssymb}
\usepackage{bm}

\newcommand{\qed}{\hfill$\blacksquare$}

\newtheorem{mydef}{Definition}
\newtheorem{mythm}{Theorem}

\newtheorem{mylemma}{Lemma}

\newtheorem{myproof}{Proof}

\newcommand{\NI}{\noindent}
\newcommand{\VS}{{\vspace*{0.1in}\NI}}

\newcommand{\reals}{\mathbb{R}}

\newcommand{\suchthat}{\;\mid\;}

\newcommand{\realsN}{\reals^N}

\newcommand{\Enorm}[1]{\| {#1} \|_2}

\newcommand{\Onenorm}[1]{\| {#1} \|_1}
\newcommand{\inverse}[1]{{#1}^{-1}}
\newcommand{\transpose}[1]{{#1}^{T}}
\newcommand{\dist}[2]{\mathrm{dist}({#1},{#2})}
\newcommand{\ssup}[2]{{#1}^{(#2)}}

\newcommand{\vv}[1]{\mathbf{#1}}
\newcommand{\vd}[1]{\dot{\mathbf{#1}}}
\newcommand{\vva}{\vv{a}}
\newcommand{\vvb}{\vv{b}}

\newcommand{\vvg}{\vv{g}}

\newcommand{\vvs}{\vv{s}}
\newcommand{\vvu}{\vv{u}}
\newcommand{\vvv}{\vv{v}}

\newcommand{\vvx}{\vv{x}}
\newcommand{\vvy}{\vv{y}}

\newcommand{\vvzero}{\vv{0}}

\newcommand{\cc}[1]{{\cal #1}}
\newcommand{\ccA}{\cc{A}}
\newcommand{\ccB}{\cc{B}}
\newcommand{\ccC}{\cc{C}}
\newcommand{\ccD}{\cc{D}}

\newcommand{\ccF}{\cc{F}}
\newcommand{\ccI}{\cc{I}}
\newcommand{\ccL}{\cc{L}}
\newcommand{\ccM}{\cc{M}}
\newcommand{\ccO}{\cc{O}}
\newcommand{\ccS}{\cc{S}}
\newcommand{\ccX}{\cc{X}}

\newcommand{\vvastar}{\vva^*}
\newcommand{\vvgstar}{\vvg^*}
\newcommand{\vvustar}{\vvu^*}
\newcommand{\vvmustar}{\boldsymbol{\mu}^*}
\newcommand{\ustar}{u^*}
\newcommand{\astar}{a^*}
\newcommand{\gstar}{g^*}
\newcommand{\mustar}{\mu^*}

\newcommand{\thres}{T}
\newcommand{\threslasso}{\thres_{\pm\lambda}}
\newcommand{\thresclasso}{\thres_\lambda}
\newcommand{\vvthres}{\vv{\thres}}

\newcommand{\vvphi}{{\boldsymbol\phi}}

\newcommand{\vvF}{\vv{F}}

\newcommand{\vdu}{\vd{u}}

\newcommand{\ustart}{\ssup{\vv{u}}{0}}
\newcommand{\vstart}{\ssup{\vv{v}}{0}}

\newcommand{\vvvhat}{\hat{\vvv}}

\newcommand{\flowlca}[2]{\vv{U}({#1},{#2})}

\newcommand{\sign}{\mathrm{sign}}

\newcommand{\fixedpts}{{\cal F}}

\newcommand{\fixedregion}{\hat{\fixedpts}}

\newcommand{\tendsto}{\rightarrow}

\newcommand{\bydef}{\stackrel{\mathrm{def}}{=}}

\begin{document}
%\title{Bare Demo of IEEEtran.cls\\ for IEEE Journals}

%%%%%============
%\input{authors}
%%%%%============

\title{Convergence of LCA Flows to (C)LASSO Solutions}

\author{
Ping Tak Peter Tang\thanks{Intel Corporation, 2200 Mission College Blvd, Santa Clara, CA 95054}
}

\maketitle

%%\input{abstract}
%%%%========

\begin{abstract}
This paper establishes several convergence results about flows of the dynamical system
LCA (Locally Competitive Algorithm) to the mixed $\ell_2$-$\ell_1$ minimization problem LASSO and
the constrained version, called CLASSO here, where the parameters are required to be non-negative.
(C)LASSO problems are closely related to various important applications including efficient coding,
image recognition and image reconstruction. That the solution of (C)LASSO can be determined
by LCA allows the former to be solved in novel ways such as through a physical realization of
analog circuits or on non-von Neumann computers. As discussed in the paper, previous works that
show convergence of LCA to LASSO are incomplete, and do not consider CLASSO. The main contributions of
this paper are a particular generalization of LaSalle's invariance principle and its application to
rigorously establish LCA's convergence to (C)LASSO.
\end{abstract}

%\begin{abstract}
%The abstract goes here.
%\end{abstract}

%\begin{IEEEkeywords}
%\end{IEEEkeywords}


% \fi
%%\IEEEpeerreviewmaketitle

%--------Start real contents
%%%%==============
%%%\input{introduction}
%%%%==============

\section{Introduction}\label{sec:introduction}

The LASSO problem (least absolute shrinkage and selection operator) is a regression problem with
regularization. It was originally formulated by Tibshirani in~\cite{Tibshirani96}. More recent works
show that solving LASSO is an important tool in various problems related to image processing,
sparse coding and compressive sensing~\cite{Fu-et-al-06}, \cite{Eldar-Kuppinger-Bolcskei-10},
\cite{Troop06}.
LASSO can be solved by traditional computational methods drawn from the
optimization algorithms~\cite{ChenDonohoSaunders01}, \cite{Kim07}, or~\cite{Boydvandenberghe04}. 
Rozell et. al.~\cite{Rozell08} formulated a dynamical system called LCA (locally competitive algorithm)
closely related to various optimization problems such as LASSO. The motivation is that instead of
using traditional numerical methods for differential equations, LCA can be solved alternatively
by running an analog circuit that is suitably configured to mimic some neural 
networks~\cite{ShaperoRozellHasler13} 
or by running a digital circuit that mimics some spiking neural networks~\cite{Kim-et-al-15}.
This approach is promising both conceptually, as it allows the use 
of non-von Neumann computing devices, and practically, as these novel devices 
can be extremely power efficient.
Consequently, firm theoretical understanding on convergence behavior of LCA to LASSO solutions is 
invaluable.

The original work~\cite{Rozell08} that formulated LCA only discussed convergence briefly. 
Furthermore, those discussions are applicable only on LCA with activation functions that are invertible.
For LCA configured for LASSO, the corresponding activation function is non-invertible, unbounded,
but not necessarily radially unbounded (see later discussions). 
Thus convergence of LCA to LASSO was not established in this original work.
Balavoine, Romberg and Rozell~\cite{Balavoine12} concurred with this assessment. After lucidly
discussing why other related works (for example, References 15 through 26 cited in~\cite{Balavoine12})
do not provide the needed convergence theory, the authors established various results, including
LCA's convergence to a LASSO solution under certain mild assumptions. 
In a later work~\cite{Balavoine13}, the same authors used the more advanced
tool of {\L}ojasiewicz inequality to not only re-establish but also strengthen their earlier convergence
results. Unfortunately, the proofs in both works have major gaps, detailed later in 
Appendix~\ref{appen:mistakes}. Thus the need for convergence guarantee remains. Moreover, in the
situation of solving LCA with the use of a spiking neural network, the LASSO parameters are 
naturally represented in terms of spiking rates, which are non-negative. Hence the corresponding LASSO
problem that LCA aims to solve is in fact a constrained version of LASSO, called CLASSO here. 
In the case of CLASSO,
the activation function in LCA differs slightly from that of LASSO. The difference is material as
this modified activation function is no longer radially unbounded, 
a property that contributes significantly to global convergence behavior of dynamical systems.

In this paper, convergence of LCA to (C)LASSO is established through a suitable generalization of the
well-known LaSalle invariance principle~\cite{LaSalle60}. In particular, the (C)LASSO objective
function value converges to its optimum along any arbitrary LCA flow (trajectory). Moreover,
when the (C)LASSO optimal solution (coefficients) is unique, the ``output'' of the LCA along any
arbitrary flow also converge to that optimal solution. 

The rest of the paper is organized as follows. Section~\ref{sec:background} reviews some standard
definitions and theories in dynamical systems as well as convex optimizations. It also states a new
generalization of the LaSalle invariance principle. Section~\ref{sec:lasso-lca} states the (C)LASSO problems
and the related LCA. In particular, theorems that relate the two problems are given.
With the set up of these two sections, Section~\ref{sec:convergence} then
describes all the convergence results of this paper. Section~\ref{sec:conclusion} 
discusses possible generalizations of the current results and an important question that is yet
unsettled. Proofs of all the technical results stated in 
Sections~\ref{sec:background} through~\ref{sec:convergence} are given in the Appendix.

%%%%==============
%\input{background}
%%%%==============

\section{Background}\label{sec:background}

This section focuses on properties related to general dynamical systems as well
as convex optimizations that will be used in the sequel. 
While most of these properties are well known and stated here for the sake of
making this paper self contained, Theorem~\ref{thm:lasalle-extension}
is new.

\subsection{Dynamical Systems}\label{ssec:dynamical-systems}
Consider the system of differential equations describing a function
$\vvu : \realsN \rightarrow \reals$:
\begin{equation}\label{eqn:DS}
\frac{d}{dt} \vvu(t) = \vvF(\vvu(t)), \quad
\hbox{or compactly} \quad
\vdu = \vvF(\vvu),\tag{DS}
\end{equation}
where $\vvF:\reals^N \rightarrow \reals^N$ is locally Lipschitz.
Standard theory of ordinary differential equations shows that the system
has a unique solution $\flowlca{t}{\ustart}$, $t \ge 0$ where
$\ustart = \flowlca{0}{\ustart}$ is a given point in $\reals^N$. 
Throughout this paper, $\flowlca{t}{\ustart}$, $t \ge 0$, is called
the flow with initial position $\ustart$.

\begin{mydef}\label{def:basics}

\NI
\begin{enumerate}
\item
A point $\vvu\in\realsN$ is called a fixed point iff $F(\vvu) = \vvzero$.
\item
A set $\ccX \subseteq \realsN$ is called positive invariant if for any
$\ustart\in\ccX$, $\flowlca{t}{\ustart}\in\ccX$ for all $t \ge 0$.
\item
Given a set $\ccM\subseteq\realsN$, a flow $\flowlca{t}{\ustart}$ is said
to converge to $\ccM$, $\flowlca{t}{\ustart}\tendsto\ccM$, if 
$$
\lim_{t\tendsto\infty}\dist{\flowlca{t}{\ustart}}{\ccM} = 0,
$$
where $\dist{\vvx}{\ccM} = \inf_{\vvy \in \ccM} \Enorm{\vvx - \vvy}$.
\item
A set $\ccD$ is called a domain of bounded flows if for each $\ustart\in\ccD$,
there is a $K > 0$ (possibly dependent on $\ustart$) such that
$\Enorm{\flowlca{t}{\ustart}} \le K$ for all $t \ge 0$.
\end{enumerate}
\end{mydef}

\NI
The following theorem is well know~\cite{LaSalle60}.

\begin{mythm}\label{thm:lasalle}
(LaSalle Invariance Principle)
Consider a dynamical system (\ref{eqn:DS}).
Let $\ccD\subseteq\realsN$ be compact and positive invariant, and
$V : \realsN \rightarrow \reals$ be a scalar function with continuous
first partial derivatives. Suppose
$$
\dot{V}(\vvu) \bydef \grad V(\vvu) \cdot \vvF(\vvu) =
\sum \frac{\partial V(\vvu)}{\partial u_n} F_n (\vvu) \le 0
$$
for all $\vvu\in\realsN$. Then for all $\ustart\in\realsN$,
$\flowlca{t}{\ustart} \tendsto \ccM$, where 
the set $\ccM$ is the largest positive invariant set contained in
$\ccS = \{\vvu \suchthat \dot{V}(\vvu) = 0\}$.
\end{mythm}

\NI
The main result of this paper relies on an extension to Theorem~\ref{thm:lasalle}.
Theorem~\ref{thm:lasalle-extension} relaxes on the compactness and smoothness
requirements on $\ccD$ and $V$, respectively.
The proof is given in Appendix~\ref{appen:invariance}.

\begin{mythm}\label{thm:lasalle-extension}
Let $\ccD\subseteq\realsN$ be a domain of bounded flows that is closed and positive
invariant. Suppose there are scalar functions $V, W:\realsN\rightarrow\reals$ such that
$V$ is continuous and $W$ is upper semicontinuous and non-positive:
$W(\vvu) \le 0$ for all $\vvu\in\realsN$. Moreover,
for every flow $\vvu(t) = \flowlca{t}{\ustart}$,
$$
\frac{d}{dt} V(\vvu(t)) = W(\vvu(t)) \quad
\hbox{a.e. in $[0, \infty)$},
$$
that is for all $t \in [0,\infty)$ except possibly for a set of measure 0. Then
for any $\ustart\in\ccD$, $\flowlca{t}{\ustart} \tendsto\ccM$, 
the largest positive invariant set inside $\ccS = \{\vvu \suchthat W(\vvu) = 0\}$.
\qed
\end{mythm}

%------------------------------------------------------------------------

\subsection{Convex Optimizations}\label{ssec:convex-optimizations}

Consider a convex function $E:\realsN\rightarrow\reals$ and the optimization problem

\begin{equation}\label{eqn:cvx-opt}
\begin{aligned}
 &\operatorname*{arg\,min}_{\vva\in\realsN} E(\vva)\\
 &\hbox{subject to}\\
 &h_i(\vva) \le 0, \quad i=1,2,\ldots,L, \; L \ge 0,
\end{aligned}\tag{CO}
\end{equation}
where each of the functions $h_i$ is affine. Note that problem~\ref{eqn:cvx-opt} 
is not the most general convex optimization problem as the only constraints are
that of inequality constraints,
which are also defined by affine functions instead of general convex functions. 
$L$ is allowed to be zero, in which case, the problem is unconstrained.

Standard theory in convex optimizations shows that for Problem~\ref{eqn:cvx-opt},
the KKT conditions are necessary and sufficient in characterizing optimal solutions
(see~\cite{Boydvandenberghe04} for details), as stated below.

\begin{mythm}\label{thm:KKT}
A point $\vvastar\in\realsN$ is an optimal solution for \ref{eqn:cvx-opt} iff there is
a $\vvmustar\in\reals^L$ such that all of the following conditions hold.
\begin{enumerate}
\item (Stationarity) 
$$\vvzero\in\partial E(\vvastar) + \sum_{i=1}^L \mu^*_i \grad h_i(\vvastar)$$
where $\partial E$ is the generalized gradient of $E$ (see \cite{Clarke90}).
\item (Complementarity)
$$
\mu^*_i\,h_i(\vvastar) = 0, \quad i=1,2,\ldots,L.
$$
\item (Feasibility)
$$
h_i(\vvastar) \le 0, \; \mu^*_i \ge 0, \quad i=1,2,\ldots,L.
$$
\end{enumerate}
If $L=0$ (the case of unconstrained optimization), these conditions reduced to simply
$\vvzero \in \partial E(\vvastar)$.
\end{mythm}

%%%%-------------------------------------------------------------------------------------------

%%%%==============
%\input{lasso-lca}
%%%%==============

\section{LASSO and LCA}\label{sec:lasso-lca}

The original LASSO problem as formulated in~\cite{Tibshirani96} is
an unconstrained convex optimization: Given $\vvs \in \reals^M$,
a matrix 
$\Phi=[\vvphi_1, \vvphi_2, \ldots,\vvphi_N] \in\reals^{M\times N}$
($\vvphi_j \in\reals^M$ is the $j$-th column) and a real number $\lambda > 0$,
solve
$$
\operatorname*{arg\,min}_{\vva\in\reals^N} \,
\frac{1}{2} \Enorm{\vvs - \Phi \vva}^2 + \lambda \Onenorm{\vva}.
$$
This paper also considers the constrained version (called CLASSO here)
where the vector of parameters $\vva$ is restricted to having only
non-negative components, denoted as $\vva \ge \vvzero$. Both versions
correspond to Problem~\ref{eqn:cvx-opt} where
$$
E(\vva) = 
\frac{1}{2} \Enorm{\vvs - \Phi \vva}^2 + \lambda \Onenorm{\vva},
$$
and $h_i(\vva) = -a_i$, $i=1,2,\ldots,L$ where either
$L=0$ (LASSO) or $L=N$ (CLASSO).

The dynamical system LCA (locally competitive algorithm) is originally
formulated to solve LASSO~\cite{Rozell08} although~\cite{Balavoine12} is a
first attempt to establish rigorously that LCA solves LASSO.
LCA can also be configured to address CLASSO. Specifically,
assuming the (C)LASSO problem is given as before but
with the $\Phi$ matrix to have columns scaled to have unit norm.
LCA is a dynamical system of the form of~\ref{eqn:DS}, $\vdu = \vvF(\vvu)$,
$\vvF:\realsN\rightarrow\realsN$ where
$$
\vvF(\vvu) = \vvb - \vvu - (\transpose{\Phi}\Phi - I)\vvthres(\vvu),
$$
$\vvb = \transpose{\Phi}\vvs$, and $\vvthres:\realsN\rightarrow\realsN$
is a nonlinear function that applies an identical scalar function $\thres:\reals\rightarrow\reals$
to each component $u_n$ of $\vvu$. The function $\thres$ is either 
$\threslasso$ for LASSO and $\thresclasso$ for CLASSO:
$$
\thresclasso(x) = 
\left\{
\begin{array}{c c}
x - \lambda & x > \lambda \\
0           & x \le \lambda 
\end{array}
\right.
$$
and $\threslasso(x) = \thresclasso(x) + \thresclasso(-x)$.

That LCA solves (C)LASSO means that an LCA flow $\vvu(t) = \flowlca{t}{\ustart}$
for some suitably (or arbitrarily) chosen initial position $\ustart$ leads
to $\vvthres(\vvu(t))$ converging to an optimal point for (C)LASSO. In practice,
this convergence usually happens in tandem with $\vvu(t)$ converging to a fixed point of LCA.
Consider a (C)LASSO problem and the corresponding LCA system.
\begin{mydef}\label{def:CFFhat}

\NI
\begin{enumerate}
\item $\ccC$ denotes the set of optimal solutions to (C)LASSO.
\item $\ccF$ denotes the set of fixed points of LCA.
\item $\fixedregion$, defined as $\fixedregion \bydef \inverse{\vvthres}(\ccC)$, is called the
      fixed region. That is a $\vvu\in\realsN$ will yield an optimal point via the $\vvthres$
      mapping iff $\vvu\in\fixedregion$.
\item Given $\vva\in\realsN$, denote the index set of nonzero and zero components by
$$
\ccA(\vva) = \{ n \suchthat a_n \neq 0 \}
$$
and
$$
\ccI(\vva) = \{ n \suchthat a_n = 0 \}.
$$
\end{enumerate}
\end{mydef}
Whenever a flow of a dynamical system converges, it necessarily converges to a fixed point. 
The next two theorems give the crucial relationship between
$\ccF$, $\fixedregion$ and $\ccC$. Their proofs are given in Appendix~\ref{appen:critical}. 

\begin{mythm}\label{thm:critical-fixedregion}
Let $\ccC$, $\ccF$ and $\fixedregion$ be as defined previously. Then the following hold:
\begin{enumerate}
\item $\vvthres(\ccF) \subseteq \ccC$.
\item For each $\vvastar\in\ccC$, there is a $\vvustar\in\ccF$ such that $\vvthres(\vvustar) = \vvastar$.
      In particular, $\ccF\subseteq\fixedregion$.
\end{enumerate}
\end{mythm}

The fixed region $\fixedregion$ plays an important role in subsequent developments of this paper. The
next theorem sets a foundation by applying the KKT characterization of Theorem~\ref{thm:KKT} to
$\fixedregion$.

\begin{mythm}\label{thm:fixedregion-KKT}
Let $\vvu\in\realsN$. Then $\vvu\in\fixedregion$ if and only if both equations below hold
$$
\begin{array}{l l}
0 = b_n - \transpose{\vvphi}_n\Phi \vva - \lambda\sign(a_n), & n\in\ccA(\vva), \\
0 = \thres\left(b_n - \transpose{\vvphi}_n\Phi\vva \right),  & n\in\ccI(\vva),
\end{array}
$$
where $\vva = \vvthres(\vvu)$.
\end{mythm}

%%%%-------------------------------------------------------------------------------------------

%%%%==============
%\input{convergence}
%%%%==============

\section{Convergence of LCA to LASSO Solutions}\label{sec:convergence}

This section establishes various convergence properties of LCA. The main tool
is the extension of LaSalle invariance principle stated in Theorem~\ref{thm:lasalle-extension}. 
In order to apply that theorem, the two scalar functions $V$ and $W$ are defined below in
Definition~\ref{def:LASSO-V-and-W} and their relevant properties are stated in 
Lemma~\ref{lemma:LASSO-V-and-W}. With $V$ and $W$ appropriately defined,
the result of Theorem~\ref{thm:lasalle-extension} says that all LCA flows converge to the
set $\ccM$, the largest positive invariant set inside the set $\ccS$ of all points at which
$W$ is zero. To further refine this result,
Theorem~\ref{thm:fixedregion-invariance} states that
this largest invariant set turns out to be the fixed region of (C)LASSO. Two
convergence results then follow easily from these foundations.
Proofs for all the results of this section are given in Appendix~\ref{appen:convergence}.

\begin{mydef}\label{def:LASSO-V-and-W}
Given a (C)LASSO problem and the corresponding LCA system, define scalar functions
$V, W:\realsN \rightarrow \reals$ as follows.
$$
V(\vvu) \bydef E(\vva) = \frac{1}{2}\Enorm{\vvs - \Phi\vva}^2 +
                              \lambda\,\Onenorm{\vva},
$$
where $\vva = \vvthres(\vvu)$, and
$$
W(\vvu) \bydef \left\{
\begin{array}{l l}
\sum_{n\in\ccA(\vva)} \frac{\partial V}{\partial u_n} \cdot F_n(\vvu) & \ccA(\vva) \neq \emptyset \\
0                                                                & \ccA(\vva) = \emptyset
\end{array}
\right.
$$
\end{mydef}

\begin{mylemma}\label{lemma:LASSO-V-and-W}
Consider $V$ and $W$ in Definition~\ref{def:LASSO-V-and-W}. The following hold.
\begin{enumerate}
\item $V$ is continuous.
\item $W$ is non-positive, that is, $W(\vvu) \le 0$ for all $\vvu\in\realsN$,
      and upper semicontinuous.
\item Given any LCA flow $\vvu(t) = \flowlca{t}{\ustart}$, $\ustart\in\realsN$,
$$
\frac{d}{dt} V(\vvu(t)) = W(\vvu(t)) \quad
\hbox{a.e. in $t\in [0,\infty)$}.
$$
\item Given any $\ustart\in\realsN$, the set
$$
\ccD = \{ \vvu\in\realsN \suchthat V(\vvu) \le V(\ustart) \}
$$
is closed, positive invariant, and a domain of bounded flows.
\end{enumerate}
\end{mylemma}

Given Definition~\ref{def:LASSO-V-and-W} and Lemma~\ref{lemma:LASSO-V-and-W}, the extension to
LaSalle's invariance principle as stated in Theorem~\ref{thm:lasalle-extension} is applicable. Thus
$\flowlca{t}{\ustart} \tendsto \ccM$ for any $\ustart\in\realsN$ where $\ccM$ is the largest
positive invariant set inside $\ccS = \{\vvu \suchthat W(\vvu) = 0 \}$. The next theorem shows that
this $\ccM$ is in fact the fixed region $\fixedregion = \inverse{\vvthres}(\ccC)$.

\begin{mythm}\label{thm:fixedregion-invariance}
Let $\ccM$ be the largest positive invariant set inside $\ccS=\{\vvu \suchthat W(\vvu) = 0 \}$.
Then $\ccM = \fixedregion$. In particular, given an arbitrary flow $\vvu(t) = \flowlca{t}{\ustart}$,
$\vvu(t) \tendsto \fixedregion$ and $\vvthres(\vvu(t)) \tendsto \ccC$.
\end{mythm}

The next theorem shows that LCA can be used to determine the optimal objective function value of (C)LASSO.
\begin{mythm}\label{thm:V-convergence}
Let $E^* = E(\vvastar)$, $\vvastar\in\ccC$, be the optimal objective function value
of (C)LASSO. Denote by
$\vvu(t) = \flowlca{t}{\ustart}$ the LCA flow at an arbitrary starting point. Then
$\lim_{t\tendsto\infty} V(\vvu(t)) = E^*$.
\end{mythm}

Because (C)LASSO is convex, the set of optimal solutions is a convex set. In particular an isolated
optimal solution exists iff the optimal solution is unique. In this case a stronger convergence
result can be established and LCA can be used to determine the optimal (C)LASSO solution.

\begin{mythm}\label{thm:flow-convergence}
Suppose (C)LASSO has a unique optimal solution $\vvastar\in\ccC$, then given any $\ustart\in\realsN$
$\vvu(t)=\flowlca{t}{\ustart}$ converges to a fixed point $\vvustar\in\ccF$: $\vvu(t)\tendsto\vvustar$.
Furthermore $\vvthres(\vvustar) = \vvastar$.
\end{mythm}

%%%%-------------------------------------------------------------------------------------------

%%%%==============
%\input{conclusion}
%%%%==============

\section{Conclusion}\label{sec:conclusion}

In summary, one can determine the optimal objective function value of (C)LASSO by
approaching the corresponding LCA's fixed region arbitrarily closely -- which any
flow $\flowlca{t}{\ustart}$ always does.
Furthermore, when the optimal (C)LASSO parameters are unique, 
$\vvthres\left(\flowlca{t}{\ustart}\right)$ will converge to that as well.
It is known (cf.~\cite{Balavoine12} and~\cite{Lu-Wang-08})
that a more general optimization problem of the form
$$
\operatorname*{arg\ min}_\vva E(\vva), \quad
E(\vva) = \frac{1}{2} \Enorm{\vvs - \Phi\vva}^2 + \sum_{n=1}^N C(a_n)
$$
is related to the LCA
$$
\dot{\vvu} = \vvb - \vvu - (\transpose{\Phi}\Phi - I)\,
\transpose{[\thres(u_1), \thres(u_2), \ldots, \thres(u_N)]}
$$
where $\vvb = \transpose{\Phi} \vvs$ and $\thres:\reals\rightarrow\reals$
is a function strictly increasing on $\{u \suchthat \thres(u) \neq 0 \}$ such that for $a \neq 0$
$C'(a) = u - a$, $a = \thres(u)$. 
The results developed here in this paper apply trivially when $C(\cdot)$ is convex. This includes for
example a LASSO-like problem called elastic net~\cite{ZouHastie05}:
$$
\operatorname*{arg\ min}_\vva E(\vva), \quad
E(\vva) = \frac{1}{2} \Enorm{\vvs - \Phi\vva}^2 + \lambda_1 \Onenorm{\vva} + \lambda_2 \Enorm{\vva}^2.
$$
The corresponding ``activation'' function $\thres:\reals\rightarrow\reals$ are
$$
\thres_{\rm ce}(u) = \left\{
\begin{array}{l l}
0  &  u \le \lambda_1 \\
\frac{u - \lambda_1}{2\lambda_2 + 1} & u > \lambda_1,
\end{array}\right.
$$
and $\thres_{\rm e}(u) = \thres_{\rm ce}(u) + \thres_{\rm ce}(-u)$.
LCA with $\thres = \thres_{\rm e}$ corresponds to the elastic net problem while
$\thres = \thres_{\rm ce}$ corresponds to adding non-negativity constraints to
the elastic net problem.

It is worth pointing our that when $C(\cdot)$ is not convex, the usefulness of
a converging LCA may be greatly diminished. In this case, while an LCA fixed
point $\vvustar$ still yield a critical point $\vvastar = \vvthres(\vvustar)$
of $E(\vva)$ in the sense that $\vvzero\in\partial E(\vvastar)$, the
last property alone is insufficient for $\vvastar$, or even $E(\vvastar)$ to be optimal.

Focusing back to (C)LASSO, 
note that while $E(\vvu(t)) \tendsto E^*$ monotonically, where $\vvu(t)$ is a LCA flow, 
Theorem~\ref{thm:flow-convergence} only guarantees $\vvu(t)\tendsto\vvustar$ and
$\vvthres(\vvu(t))\tendsto\vvastar$ should the optimal $\vvastar$ is unique. 
While~\cite{Balavoine13} states that both $\vvu(t)$ and $\vvthres(\vvu(t))$ converge to some
$\vvustar$ and $\vvastar$, the proof's gap outlined in Appendix~\ref{appen:mistakes} renders
the convergence claims, however welcome, unsupported. Settling the convergence question, affirmatively
or otherwise, is a natural next step to the results given in this paper.

%%%%-------------------------------------------------------------------------------------------

%%\appendices
\appendix

\renewcommand{\theequation}{A-\arabic{equation}}
\setcounter{equation}{0}

\section{Proof Of Theorem~\ref{thm:lasalle-extension}}\label{appen:invariance}

Let $\ustart\in\ccD$ be an arbitrary starting point and $\ccL$ be the set of limit points
of the flow $\flowlca{t}{\ustart}$. That is $\ccL$ is the limit points of the set
$\{\flowlca{t}{\ustart} \suchthat t \ge 0\}$. The first lemma establishes some basic
topological properties of $\ccL$.

\begin{mylemma}\label{lemma:topological-properties-of-limit-set}
The limit set $\ccL$ has the following properties. 
\begin{enumerate}
\item $\ccL$ a compact subset of $\ccD$. 
\item $\ccL$ is positive invariant.
\item $\flowlca{t}{\ustart} \tendsto \ccL$ as $t \tendsto \infty$.
\end{enumerate}
\end{mylemma}
\begin{myproof}
Proof of (1):
Since $\ccD$ is a domain of bounded flows, there is a $K$ such that $\Enorm{\flowlca{t}{\ustart}}\le K$
for all $t\ge 0$. By Bolzano-Weiestrass theorem, $\ccL$ is nonempty where each of its elements must be
bounded by $K$. Because $\ccD$ is closed, we must have $\ccL \subseteq \ccD$. To show that $\ccL$ is
compact, it suffices to show that it is closed. Consider
any sequence $\{\vvu_k\}$, $\vvu_k \in \ccL$ for all $k$, that is convergent to a certain
limit $\hat{\vvu} \in \reals^N$. If $\hat{\vvu}$ must necessarily belong to $\ccL$ as well, then $\ccL$
is closed. Since $\vvu_1 \in \ccL$, there is a time $t_1 \ge 1$ such that
$\Enorm{\flowlca{t_1}{\ustart} - \vvu_1} < 1$. Having obtained $t_1$, there must be a
$t_2 > \max\{t_1,2\}$ such that $\Enorm{\flowlca{t_2}{\ustart} - \vvu_2} < 1/2$. Continuing this process,
one obtains a sequences $0 \le t_1 < t_2 < \cdots$, $t_k \tendsto\infty$, that satisfies
$\Enorm{\flowlca{t_k}{\ustart} - \vvu_k} < 1/k$.
Now for any $\epsilon > 0$, pick $K_1$ large enough so that $\Enorm{\vvu_k - \hat{\vvu}} < \epsilon/2$
for all $k \ge K_1$. Pick integer $K_2$ so that $K_2 \ge 2/\epsilon$, that is, $1/k < \epsilon/2$
for all $k \ge K_2$. Let $K = \max\{K_1,K_2\}$. Then for all $k \ge K$,
$$
\Enorm{\flowlca{t_k}{\ustart} - \hat{\vvu}}
\le 
\Enorm{\flowlca{t_k}{\ustart} - \vvu_k} + 
\Enorm{\vvu_k - \hat{\vvu}} < \epsilon.
$$
Thus $\flowlca{t_k}{\ustart} \tendsto \hat{\vvu}$ and $\hat{\vvu}$ must also belong to $\ccL$,
showing that $\ccL$ is closed. 

\VS
Proof of (2): Let $\vstart \in \ccL$ be an arbitrary starting point, and let
$\tau > 0$ be an arbitrary time. Let $\vvvhat = \flowlca{\tau}{\vstart}$. It
suffices to show $\vvvhat \in \ccL$. Since $\vstart \in \ccL$, there is a time sequence
$\{t_k\}$ such that $\vvu_k \tendsto \vstart$, where
$\vvu_k = \flowlca{t_k}{\ustart}$. 
Consider the fact
$$
\flowlca{t_k+\tau}{\ustart} =
\flowlca{\tau}{\flowlca{t_k}{\ustart}} =
\flowlca{\tau}{\vvu_k}.
$$
Now, define $\hat{t}_k = t_k + \tau$. 
Since $\vvu_k \tendsto \vstart$ and $\flowlca{\tau}{\vvu_k}$ is continuous in its
second argument, 
$$
\flowlca{\hat{t}_k}{\ustart} \tendsto \flowlca{\tau}{\vstart} = \hat{\vvv}.
$$
This shows $\hat{\vvv} \in \ccL$, establishing positive invariance of $\ccL$.

\VS
Proof of (3):
Assume the contrary.
This means that there is a $\epsilon > 0$ and a time sequence $\{t_k\}$ such that
for all $k=1,2,\ldots$,
$\dist{\flowlca{t_k}{\ustart}}{\ccL} \ge \epsilon$.
But $\{\flowlca{t_k}{\ustart}\}$ is bounded  and must have a convergent subsequence
$\flowlca{\hat{t}_k}{\ustart} \tendsto \hat{\vvu}$. Therefore, on the one hand,
$\hat{\vvu} \in \ccL$, that is, $\dist{\hat{\vvu}}{\ccL} = 0$, but on the other
hand $\dist{\hat{\vvu}}{\ccL} \ge \epsilon$ because 
$\dist{\flowlca{\hat{t}_k}{\ustart}}{\ccL} \ge \epsilon$ for all $k$.
This contradiction shows that indeed 
$\flowlca{t}{\ustart}$ converges to $\ccL$.
\qed
\end{myproof}

The second lemma concerns properties of the scalar functions $V$ and $W$.
\begin{mylemma}\label{lemma:V-and-W}
The function $V(\vvu)$ takes on a constant value on $\ccL$ and $W(\vvu) = 0$
for all $\vvu\in\ccL$.
\end{mylemma}
\begin{myproof}
Observe that for any $0\le t_1 \le t_2$,
\begin{eqnarray*}
\lefteqn{V(\flowlca{t_2}{\ustart})} \\
& = & V(\flowlca{t_1}{\ustart}) + \int_{t_1}^{t_2} W(\flowlca{s}{\ustart})\,ds \\
& \le & V(\flowlca{t_1}{\ustart}).
\end{eqnarray*}
Thus $V(\flowlca{t}{\ustart})$ is non-increasing in $t$. 
Now assume that $V(\vvu)$ is not constant on $\ccL$:
Let $\ssup{\vvu}{1}$ and $\ssup{\vvu}{2} \in \ccL$ be such that
$$
V(\ssup{\vvu}{1}) = \ssup{V}{1} < \ssup{V}{1} + \delta
= \ssup{V}{2} = 
V(\ssup{\vvu}{2})
$$
for some $\delta > 0$. $\ssup{\vvu}{1} \in \ccL$ implies there is a time $t_1$
large enough that the proximity of $\flowlca{t_1}{\ustart}$ to $\ssup{\vvu}{1}$
implies 
$$
\left| V(\flowlca{t_1}{\ustart}) - \ssup{V}{1} \right| < \delta/3.
$$
Based on $t_1$, pick $t_2 > t_1$ such that the proximity of
$\flowlca{t_2}{\ustart}$ to $\ssup{\vvu}{2}$ implies 
$$
\left| V(\flowlca{t_2}{\ustart}) - \ssup{V}{2} \right| < \delta/3.
$$
Thus we have $t_1 < t_2$ while
$$
V(\flowlca{t_1}{\ustart}) < V(\flowlca{t_2}{\ustart}) - \delta/3,
$$
contradiction the fact that $V(\flowlca{t_1}{\ustart})$ is
non-increasing in time. This establishes the first assertion.

\VS
To prove the second assertion, assume the contrary that $W(\ssup{\vvu}{1}) \neq 0$
for some $\ssup{\vvu}{1}\in\ccL$. Therefore 
$W(\ssup{\vvu}{1}) = -\epsilon$ for
some $\epsilon > 0$ as $W$ is assumed to be non-positive on $\realsN$.
Because $W$ is also assumed upper semicontinuous, there exist
positive $\epsilon$ and  $\delta_u0$ such that 
$$
W(\vvu) \le -\epsilon/2
\quad \hbox{for all
$\Enorm{\vvu - \ssup{\vvu}{1}} < \delta_u$.}
$$
By continuity of flows, there is a $\delta_t > 0$ such that
$$
\Enorm{\flowlca{t}{\ssup{\vvu}{1}} - \ssup{\vvu}{1}} < \delta_u
\quad \hbox{for all $0 \le t \le \delta_t$}.
$$
Hence $W(\flowlca{t}{\ssup{\vvu}{1}}) \le -\epsilon/2$
for all $0 \le t \le \delta_t$. Let
$\ssup{\vvu}{2} = \flowlca{\delta_t}{\ssup{\vvu}{1}}$. 
Note that $\ssup{\vvu}{2}\in\ccL$ because the $\ccL$ is positive invariant. Thus
as established previously that $V$ takes on constant value
on $\ccL$, $V(\ssup{\vvu}{2}) = V(\ssup{\vvu}{1})$.
On the other hand, however,
\begin{eqnarray*}
\lefteqn{
V(\ssup{\vvu}{2}) =
V(\ssup{\vvu}{1}) + \int_0^{\delta_t} W(\flowlca{s}{\ssup{\vvu}{1}})\,ds} \\
 & \le & V(\ssup{\vvu}{1}) - \frac{\epsilon}{2}\delta_t < V(\ssup{\vvu}{1}).
\end{eqnarray*}
This is a contradiction and thus in fact $W(\vvu) = 0$ for all $\vvu\in\ccL$.
\qed
\end{myproof}

%%%%-------------------------------------------------------------------------------------------

%\input{appen-critical}

\section{Proof Of Theorems in Section~\ref{sec:lasso-lca}}\label{appen:critical}

\begin{myproof}(For Theorem~\ref{thm:critical-fixedregion})

To establish the first part of Theorem~\ref{thm:critical-fixedregion}, $\vvthres(\ccF) \subseteq \ccC$,
it suffices to show that $\vvastar \bydef \vvthres(\vvustar) \in \ccC$ for each $\vvustar\in\ccF$.
Given a $\vvustar\in\ccF$ and defining $\vvastar = \vvthres(\vvustar)$, let the $N$ components of these
vectors be $\ustar_j, \astar_j$, $j=1,2,\ldots,N$.
Define the vector $\vvmustar, \vvgstar \in \realsN$ componentwise via
$$
\mustar_n = 
\left\{\begin{array}{c l}
0 & n\in\ccA(\vvastar) \\
-\threslasso(\ustar_n) & n\in\ccI(\vvastar)
\end{array}\right.
$$
and
$$
\gstar_n =
\left\{\begin{array}{l l}
-b_n + \transpose{\vvphi}_n\Phi\vvastar + \lambda\,\sign(\astar_n) - \mustar_n & n\in\ccA(\vvastar) \\
-b_n + \transpose{\vvphi}_n\Phi\vvastar + 
\max\{\ustar_n, -\lambda\} - \mustar_n & n\in\ccI(\vvastar) \\
\end{array}\right. 
$$
Note that $\vvmustar \ge \vvzero$ and that $\vvmustar = \vvzero$ for LASSO. 
Because 
$$
\lambda\,\sign(\astar_n) = \ustar_n - \astar_n, \quad \hbox{for $n\in\ccA(\vvastar)$, and}
$$
$$
\max\{\ustar_n,-\lambda\}-\mustar_n = \ustar_n \quad \hbox{for $n\in\ccI(\vvastar)$,}
$$
$\vvgstar = -\vvF(\vvustar) = \vvzero$. Observe now that 
$\vvgstar \in \partial E(\vvastar) - \vvmustar$. To see this, first
examine $\partial E(\vvastar)$. 
\begin{eqnarray*}
\partial E(\vva) 
 & = &
 \grad\left( \frac{1}{2}\Enorm{\vvs - \Phi \vva}^2 \right) + \lambda \partial \Onenorm{\vva} \\
 & = & -\vvb + \transpose{\Phi}\Phi \vva + \lambda\,\partial \Onenorm{\vva},
\end{eqnarray*}
where the $n$-th component of $\partial\Onenorm{\vva}$ is either $\{\sign(a_n)\}$ if
$|a_n| > 0$ or $[-1,1]$ for $a_n = 0$. Clearly then
$$
\gstar_n \in \left\{
\begin{array}{l l}
\{-b_n + \transpose{\vvphi}_n\Phi\vvastar + \lambda\,\sign(\astar_n) - \mustar_n \} & n\in\ccA(\vvastar) \\
-b_n + \transpose{\vvphi}_n\Phi\vvastar + [-\lambda,\lambda] - \mustar_n          & n\in\ccI(\vvastar)
\end{array}
\right.,
$$
that is, $\vvgstar \in\partial E(\vvastar) - \vvmustar$, implying
$\vvzero \in\partial E(\vvastar) - \vvmustar$.
Stationarity of KKT (Theorem~\ref{thm:KKT}) is satisfied for (C)LASSO. In the case of
CLASSO, the complementarity and feasibility conditions are also satisfied as 
$-\mustar_n\,\astar_n = 0$ for all $n$ and $\vvmustar\ge\vvzero$. Consequently,
$\vvastar\in\ccC$ and the first part of Theorem~\ref{thm:critical-fixedregion} is proved.

To prove the second part, let $\vvastar\in\ccC$. 
In the case of CLASSO, because $\vvastar\in\ccC$, there is a 
$\vvmustar\in\realsN$ that satisfies the three KKT conditions. In the case of LASSO,
simply define $\vvmustar$ to be the zero vector in $\realsN$. By stationarity of KKT,
\begin{equation}\label{eqn:zero-inside}
\vvzero \in
\left\{
\begin{array}{l l}
\{-b_n + \transpose{\vvphi}_n\Phi\vvastar + \lambda\,\sign(\astar_n) \}           & n\in\ccA(\vvastar) \\
-b_n + \transpose{\vvphi}_n\Phi\vvastar + [-\lambda,\lambda] - \mustar_n          & n\in\ccI(\vvastar)
\end{array}
\right.
\end{equation}
Define $\vvustar\in\realsN$ as follows. For $n\in\ccA(\vvastar)$, 
define $\ustar_n = \inverse{\thres}(\astar_n)$ 
as $\astar_n \neq 0$. Note that in this case $\lambda\,\sign(\astar_n) = \ustar_n - \astar_n$.
For $n\in\ccI(\vvastar)$, Equation~\ref{eqn:zero-inside} shows that there is $\alpha_n\in[-1,1]$
such that
$$
0 = -b_n + \transpose{\vvphi}_n\Phi\vvastar + (\alpha_n - \mustar_n).
$$
Define $\ustar_n = \alpha_n - \mustar_n$ and thus $\thres(\ustar_n) = 0$. Consequently,
$\vvthres(\vvustar) = \vvastar$ and $\vvF(\vvustar) = \vvzero$ and the second part of this
theorem is proved.\qed
\end{myproof}

\begin{myproof}(For Theorem~\ref{thm:fixedregion-KKT})

By definition of $\fixedregion$, $\vvu\in\fixedregion$ iff
$\vva=\vvthres(\vvu)\in\ccC$, iff $\vva$ satisfies KKT, iff $\vva$ is feasible
and the following two conditions hold:
$$
\begin{array}{l l l l}
0 & =   & -b_n + \transpose{\vvphi}_n\Phi \vva + \lambda\,\sign(a_n)          & n \in \ccA(\vva), \\
0 & \in & -b_n + \transpose{\vvphi}_n\Phi \vva + [-\lambda,\lambda] - \mu_n   & n \in \ccI(\vva), \\ 
\end{array}
$$
where $\mu_n = 0$ for LASSO and $\mu_n \ge 0$ for CLASSO. But the second condition is equivalent to
$b_n - \transpose{\vvphi}_n\Phi \vva \in [-\lambda,\lambda] - \mu_n$, which is equivalent to
$$
\thres\left(b_n - \transpose{\vvphi}_n\Phi \vva\right) = 0
$$
where $\thres = \threslasso$ for LASSO and $\thres =\thresclasso$ for CLASSO. This completes the
proof.\qed
\end{myproof}

%%%%-------------------------------------------------------------------------------------------

%\input{appen-convergence}

\section{Proofs for Section~\ref{sec:convergence}}\label{appen:convergence}

\begin{myproof}(For Lemma~\ref{lemma:LASSO-V-and-W}.)

Proof of (1):
That $V$ is continuous is obvious because $\Enorm{\cdot}$, $\Onenorm{\cdot}$, and
$\vvthres(\cdot)$ are all continuous.

Proof of (2): Now define for $n=1,2,\ldots,N$, 
$\chi_n(\vvu) = 1$ if $|\thres(u_n)| > 0$ and 0 otherwise. Note that
$W(\vvu) = \sum_{n=1}^N -F_n^2(\vvu)\cdot \chi_n(\vvu)$. Thus
$W(\vvu)\le 0$ for all $\vvu\in\realsN$. Observe that $\chi_n$ is lower
semicontinuous, and thus $F_n^2(\vvu) \chi_n(\vvu)$ is lower
semicontinuous, and its negation upper semicontinuous. Being a
sum of upper semicontinuous function, $W(\vvu)$ is upper semicontinuous as
well.

Proof of (3):
Let $\ccO_n = \{ t \suchthat |\thres(u_n(t))| > 0 \}$. It is open
and thus can be expressed as a countable union of disjoint open intervals
$\ccO_n = \cup_{j=1}^\infty (\alpha_j, \beta_j)$. Let
$\ccB_n = \{\alpha_j\}_1^\infty \cup \{\beta_j\}_1^\infty$.
For any $t_0 \in \ccB_n^c$ (complement of $\ccB_n$) there is an open neighborhood
$(t_0-\epsilon, t_0+\epsilon)$ in which $\frac{d}{dt}T(u_n(t))$ exists. This
derivative equals 1 if $|\thres(u_n(t))| > 0$; it equals 0 if $t_0$ is in the
interior of the set $\{t \suchthat \thres(u_n(t))=0 \}$.
The set $\ccO = \left( \cup_{n=1}^N \ccB_n \right)^c$ is open and consists of
the set $[0,\infty)$ except possibly for a set of measure 0. For each $t\in\ccO$,
$V(\vvu(t))$ is differentiable and $\frac{d}{dt} V(\vvu(t)) = W(\vvu(t))$.

Proof of (4):
Let $\ustart\in\realsN$ be chosen arbitrarily and define 
$\ccD=\{\vvu \suchthat V(\vvu)\le V(\ustart)\}$. Denote
$V(\ustart)$ by $\ssup{\eta}{0}$. $\ccD$ is closed
because it is $\inverse{V}([0,\ssup{\eta}{0}])$ and $V$ is continuous.
Since $\dot{V}(\vvu(t)) = W(\vvu(t))$ a.e. and $W(\vvu) \le 0$ always,
$V(\vvu(t))$ is non-increasing in $t$ for any flow $\vvu(t)$. Therefore
given any $\vstart\in\ccD$ and $t \ge 0$, we have
$$
V(\flowlca{t}{\vstart}) \le V(\flowlca{0}{\vstart}) \le \ssup{\eta}{0},
$$
implying that $V(\flowlca{t}{\vstart})\in\ccD$ and $\ccD$'s
positive invariance. Finally, for any $\vstart\in\ccD$, let
$\vvu(t) = \flowlca{t}{\vstart}$ be the flow with initial position 
$\vstart$ inside $\ccD$. Because $V(\vvu(t)) \le \ssup{\eta}{0}$ for
all $t \ge 0$ and $V(\vvu) \ge \lambda \Onenorm{\vvthres(\vvu)}$, there
is a constant $K_1$ such that 
$\Enorm{\vvb + (\transpose{\Phi}\Phi - I)\vvthres(\vvu(t))} \le K_1$
for all $t\ge 0$.
Since $\dot{\vvu} = -\vvu + (\vvb + (\transpose{\Phi}\Phi - I)\vvthres(\vvu))$,
$\vvu(t)$ satisfies
$$
\vvu(t) = \vstart + e^{-t}
          \int_0^t e^s \left(\vvb + (\transpose{\Phi}\Phi-I)\vvthres(\vvu(s))\right)\,ds.
$$
Therefore there is $K\ge 0$ such that $\Enorm{\vvu(t)} \le K$ for all $t\ge 0$,
establishing that $\ccD$ is a domain of bounded flows.
\end{myproof}

The following Lemma is instrumental to proving Theorem~\ref{thm:fixedregion-invariance}.
\begin{mylemma}\label{lemma:trajectory-V}
Let $\hat{\vvu}\in\realsN$ and $\hat{\vva}=\vvthres(\hat{\vvu})$ such that
$$
0 = -b_n + \transpose{\vvphi}_n\Phi\hat{\vva} + \lambda\,\sign(\hat{a}_n), \quad
\hbox{for all $n\in\ccA(\hat{\vva})$}.
$$
Define the function $\vvv(t)$ by
$$
v_n(t) = \left\{
\begin{array}{l l}
\hat{u}_n \quad \hbox{for all $t\ge 0$}  & n \in \ccA(\hat{\vva}) \\
\alpha_n
-e^{-t}( \alpha_n - \hat{u}_n )          & n \in \ccI(\hat{\vva})
\end{array}
\right.
$$
where
$\alpha_n\bydef b_n-\transpose{\vvphi}_n\Phi\hat{\vva}$.
If there exists a $\tau>0$ such that 
$\thres(v_n(t)) = 0$ for all $t\in[0,\tau]$ and
for all $n\in\ccI(\hat{\vva})$, then
$\vvv(t) = \flowlca{t}{\hat{\vvu}}$ for all $t\in[0,\tau]$.
\end{mylemma}
\begin{myproof}(For Lemma~\ref{lemma:trajectory-V})

By construction $\vvv(0) = \hat{\vvu}$. Let $\tau>0$
be as described. Then by definition of $\vvv(t)$, we
have $\vvthres(\vvv(t)) = \hat{\vva}$ for
$t\in[0,\tau]$. For $n\in\ccI(\hat{\vva})$ and
$t\in[0,\tau]$,
\begin{eqnarray*}
0 & = & (b_n-\transpose{\vvphi}_n\Phi\hat{\vva}) - v_n(t) \\
  & = & b_n - v_n(t) - \transpose{\vvphi}_n\Phi\vvthres(\vvv(t)) +
        \thres(v_n(t)) \\
  & = & F_n(\vvv(t)).
\end{eqnarray*}
For $n\in\ccA(\hat{\vva})$ and all $t\ge 0$,
\begin{eqnarray*}
dot{v}_n(t) & = & 0 \\
& = & b_n - \transpose{\vvphi}_n \Phi\hat{\vva} - \lambda\,\sign(\hat{a}_n) \\
& = & b_n - \transpose{\vvphi}_n \Phi\hat{\vva} - (\hat{u}_n - \hat{a}_n) \\
& = & F_n(\vvv(t)).
\end{eqnarray*}
Hence $\vvv(t) = \vvF(\vvv(t))$ for $t\in [0,\tau]$ and
$\vvv(0) = \hat{\vvu}$, that is,
$\vvv(t) = \flowlca{t}{\hat{\vvu}}$ for $t\in[0,\tau]$.\qed
\end{myproof}

\begin{myproof}(For Theorem~\ref{thm:fixedregion-invariance})
$\fixedregion = \ccM$ is equivalent to 
$\fixedregion\subseteq\ccM$ and $\ccM\subseteq\fixedregion$. To
prove $\fixedregion\subseteq\ccM$, consider $\hat{\vvu}\in\fixedregion$.
Denote $\vvthres(\hat{\vvu})$ by $\hat{\vva}$.
By Theorem~\ref{thm:fixedregion-KKT},
$$
0 = \left\{
\begin{array}{l l}
-b_n + \transpose{\vvphi}_n\Phi \hat{\vva} + \lambda\,\sign(\hat{a}_n) & n \in \ccA(\hat{\vva}) \\
\thres\left(b_n + \transpose{\vvphi}_n\Phi \hat{\vva}\right)           & n \in \ccI(\hat{\vva}).
\end{array}
\right.
$$
Define $\vvv(t)$ as in Lemma~\ref{lemma:trajectory-V}:
\begin{equation}\label{eqn:trajectory-V}
v_n(t) = \left\{
\begin{array}{l l}
\hat{u}_n                               &  n \in \ccA(\hat{\vva}) \\
\alpha_n - e^{-t}(\alpha_n - \hat{u}_n) &  n \in \ccI(\hat{\vva})
\end{array}
\right.
\end{equation}
where $\alpha_n \bydef b_n-\transpose{\vvphi}_n\Phi\hat{\vva}$. 
Because $\thres(\alpha_n) = 0$ for all $n\in\ccI(\hat{\vva})$
and $\inverse{\thres}(\{0\})$ is connected (either $[-\lambda,\lambda]$ or
$(-\infty,\lambda]$), $\thres(v_n(t)) = 0$ for all $t\ge 0$ and $n\in\ccI(\hat{\vva})$.
By Lemma~\ref{lemma:trajectory-V}, $\vvv(t) = \flowlca{t}{\hat{\vvu}}$ for all $t\ge 0$.
Observe that
$$
W(\vvv(t)) = \sum_{n\in\ccA(\vvthres(\vvv(t)))} -F_n^2(\vvv(t)) \ge
             \sum_{n\in\ccA(\hat{\vva})} -F_n^2(\vvv(t))  = 0,
$$
which implies $W(\vvv(t)) = 0$ for all $t\ge 0$ as $W(\vvu) \le 0$ for all
$\vvu\in\realsN$. Consequently, $\hat{\vvu}\in\ccM$ and $\fixedregion\subseteq\ccM$
is established.

Now consider a $\hat{\vvu}\in\ccM$ and $\hat{\vva} = \vvthres(\hat{\vvu})$. Thus
$W(\hat{\vvu}) = 0$ and we must have
$$
0 = F_n(\hat{\vvu}) = b_n - \hat{u}_n - \transpose{\vvphi}_n\Phi\hat{\vva} + \hat{a}_n
$$
for all $n\in\ccA(\hat{\vva})$. This is equivalent to
$$
0 = -b_n + \transpose{\vvphi}_n\Phi\hat{\vva} + \lambda\,\sign(\hat{a}_n)
$$
for $n\in\ccI(\hat{\vva})$. Define $\vvv(t)$ as in (\ref{eqn:trajectory-V}).
Because $\thres(v_n(0)) = \thres(\hat{u}_n) = 0$ for all $n\in\ccI(\hat{\vva})$,
we have $\tau > 0$ where
$$
\tau \bydef \inf \{ t \suchthat  |\thres(v_n(t))| > 0, n \in\ccI(\hat{\vva}) \}.
$$
By Lemma~\ref{lemma:trajectory-V}, $\vvv(t) = \vvu(t) \bydef \flowlca{t}{\hat{\vvu}}$
for all $0 \le t \le \tau$. If $\tau$ is finite, then there is an $n\in\ccI(\hat{\vva})$
such that $|v_n(\tau)| = \lambda$ and $\dot{v}_n(t) v_n(t) > 0$. Thus we also have
$|u_n(\tau)| = \lambda$ and $\dot{u}_n(t) u_n(t) > 0$. By continuity of $u_n(t)$
and $\dot{u}_n(t)$, there is a $\delta > 0$ such that
$$
|\thres(u_n(\tau+\delta))| > 0 \quad {\rm and} \quad
|\dot{u}_n(\tau+\delta)| > 0.
$$
Consequently, $W(u(\tau+\delta)) < 0$, contradicting the assumption that $\hat{u}\in\ccM$
which is a positive invariant set where $W(\vvu) = 0$ for all $\vvu\in\ccM$.
Consequently $\tau$ must in fact be infinite. This means that $\thres(v_n(t)) = 0$
for all $t\ge 0$ and all $n\in\ccI(\hat{\vva})$. Consequently,
$\thres(\alpha_n) = 0$ for all $n\in\ccI(\hat{\vva})$, where
$\alpha_n = b_n - \transpose{\vvphi}_n\Phi\hat{\vva}$. Note also that
$W(\hat{\vvu}) = 0$ implies for all $n\in\ccA(\hat{\vva})$,
\begin{eqnarray*}
0 & = & F_n(\hat{\vvu}) \\
  & = & b_n - \hat{u}_n - \transpose{\vvphi}_n\Phi\hat{\vva} + \hat{a}_n \\
  & = & -b_n + \transpose{\vvphi}_n\Phi\hat{\vva} + \lambda,\sign(\hat{a}_n).
\end{eqnarray*}
Using Theorem~\ref{thm:fixedregion-KKT}, one concludes that $\hat{\vvu}\in\fixedregion$.
Thus $\ccM\subseteq\fixedregion$. Together with the previously
established fact $\fixedregion\subseteq\ccM$, $\fixedregion = \ccM$ and the proof
of this theorem is complete.

Now that $\fixedregion = \ccM$, the convergence of an arbitrary flow $\vvu(t)$ 
to $\fixedregion$, $\vvu(t)\tendsto\fixedregion$ follows immediately from 
Theorem~\ref{thm:lasalle-extension}. Finally, because
$\vvthres(\cdot)$ is uniformly continuous, we have 
$\vvthres(\vvu(t)) \tendsto \vvthres(\fixedregion) = \ccC$.
\end{myproof}

\begin{myproof}(For Theorem \ref{thm:V-convergence})
This theorem follows easily from Theorem~\ref{thm:fixedregion-invariance}. Given
an arbitrary flow $\vvu(t) = \flowlca{t}{\ustart}$, $\ustart\in\realsN$, 
Theorem~\ref{thm:fixedregion-invariance} states that $\vvthres(\vvu(t)) \tendsto \ccC$.
Thus for all $t \ge 0$, there is $\ssup{\vva}{t}\in\ccC$ such that
$$
\epsilon_t = \Enorm{\vvthres(\vvu(t)) - \ssup{\vva}{t}} \tendsto 0
\quad \hbox{as $t \tendsto \infty$}.
$$
Note that 
$$
E(\vva) = \frac{1}{2} \Enorm{ \vvs - \Phi \vva }^2 + \lambda\,\Onenorm{\vva}
$$
is clearly Lipschitz in $\vva$. Thus, there is a $K > 0$ such that
$$
\left| V(\vvu(t)) - E^* \right| =
\left| E(\vvthres(\vvu(t))) -  E(\ssup{\vva}{t}) \right|
\le K \epsilon_t \tendsto 0.
$$
This completes the proof.
\qed
\end{myproof}

\begin{myproof}(For Theorem \ref{thm:flow-convergence})
Let $\vvastar$ be the unique (C)LASSO solution. By Theorem~\ref{thm:fixedregion-invariance} then,
any flow $\vvu(t)=\flowlca{t}{\ustart}$, $\ustart\in\realsN$, has the convergence property
$\vvthres(\vvu(t)) \tendsto \ccC$, the latter of which is the singleton  $\{\vvastar\}$ by 
assumption. Thus $\vvthres(\vvu(t)) \tendsto \vvastar$. Consequently, the LCA
$$
\dot{\vvu}(t) = \vvb - \vvu(t) - (\transpose{\Phi}\Phi - I)\,\vvthres(\vvu(t))
$$
must necessarily converge to the fixed point
$$
\vvustar = \vvb - (\transpose{\Phi}\Phi - I)\,\vvastar.
$$
This completes the proof.
\qed
\end{myproof}

%%%%-------------------------------------------------------------------------------------------

%\input{appen-mistakes}

\section{Technical Gaps in \cite{Balavoine12} and \cite{Balavoine13}}\label{appen:mistakes}

As alluded to before, the paper~\cite{Balavoine12} pointed out the need for theoretical results
on the convergence of LCA and proceeded to provide them.
The formulation involves a thresholding function denoted as $\vvthres$.
While the components of $\vvthres$ in~\cite{Balavoine12} is more general than $\threslasso$ here,
The difference is irrelevant for the purpose of discussion here.
In the notation of~\cite{Balavoine12}, the LCA is
$$
\dot{\vvu}(t) = \vvb - \vvu(t) - (\transpose{\Phi}\Phi - I)\vva(t), \quad
\vva(t) = \vvthres(\vvu(t))
$$
and $V(\vvx)$ is the corresponding ``energy'' function
$$
V(\vvx) = \frac{1}{2} \Enorm{\vvs - \Phi \vvx} + \lambda \sum_{m=1}^N C(x_m).
$$
Again, for the purpose of the discussion here, $C(x)$ can be considered to be the
function $|x|$. Note that this $V$ function is different from the one 
Definition~\ref{def:LASSO-V-and-W}.

The main results are found in Theorem 1 which states two convergence behavior of 
$\vva(t) = \vvthres(\vvu(t))$, $\vvu(t)$ being a LCA flow, when the 
LASSO critical points are isolated.
The proof (presented in the appendix of the paper) states that $V(\vva(t)) \tendsto V^*$
for some $V^*$ because $V(\vva(t))$ is non-increasing and bounded below. This is a sound
statement. The paper uses this fact to infer that (a) $\dot{V}(\vva(t)) \tendsto 0$,
(b) $\Enorm{\dot{\vva}(t)} \tendsto 0$, and (c) $\vva(t)\tendsto {\cal X}$ where
${\cal X}$ is the set $\{ \vva | \dot{\vva}(t) = 0 \}$.
The arguments in (b) and (c) depend crucially on (a). The argument supporting (a)
was that $V(\vva(t))$ being non-increasing and bounded below. (To quote verbatim:
``...$V(\vva(t))$ is nonincreasing for all $t\ge 0$. Since $V(\vva(t))$ is continuous,
bounded below by zero, and nonincreasing, $V(\vva(t))$ converges to a constant value
$V^*$, and its time derivative $\dot{V}(\vva(t))$ tends to zero as $t\tendsto\infty$.'')
This argument is invalid: Consider a continuous function $V(t)$, $V(0) = 0$ such that
$\dot{V}(t) = 0$ for all $[0,\infty)$ except taking on the value $-1$ on the
intervals $(k, k+1/k^2)$ for $k = 1, 2, \dots$. While $V(t)$ is nonincreasing
and bounded below (by $-\sum 1/k^2$), $\dot{V}(t)$ cannot converge to 0. 

As the technical arguments in the sequel of (a) depend crucially on it, the convergence proof
in~\cite{Balavoine12} is invalid.

In~\cite{Balavoine13} the authors of~\cite{Balavoine12} use a more advanced mathematical
tool to prove a stronger result. In addition to the original convergence result 
in~\cite{Balavoine12}, this work proves that the LCA flow $\vvu(t)$ (not just 
$\vva(t)$) converges to a $\vvu^*\in\reals^N$, regardless the critical points of $V(\vva)$
or LCA fixed points are isolated or not. This result is stronger than the one established in
this current paper. One crucial relationship that allows this advanced mathematical
tool to be used successfully is in Section IV: 
\begin{equation}\label{eqn:Ba13_1}
\begin{aligned}
\Enorm{\dot{\vva}(t)} 
  & =    \Enorm{\dot{\vva}_\Gamma(t)} \ge \beta \Enorm{\dot{\vvu}_\Gamma(t)} \\
  & \ge  \beta\,m(\partial V(\vva_\Gamma(t))) =
          \beta\,m(\partial V(\vva(t))).
\end{aligned}
\end{equation}
Here $\Gamma$ is the set of active nodes (which is $\ccA(\vva(t))$ in the current paper),
$\partial V(\vva)$ is the generalized gradient at a specific point $\vva$, 
which is in general a set of vectors, and $m( \partial V(\vva) )$ is 
$$
m(\partial V(\vva)) \bydef \inf \{ \Enorm{\vvg} \suchthat \vvg \in \partial V(\vva) \}.
$$
The $\alpha$ and $\beta$ are positive constants related to the specific $\vvthres$ being used.
As long as they are positive, their specific values do
not matter as far as the discussion below is concerned. The technical relationship 
in Equation~\ref{eqn:Ba13_1} above is
crucial as it allows $\dot{V}(\vva(t))$ be bounded in terms of 
$\Enorm{\vva(t)}\,m(\partial V(\vva(t)))$. The {L}ojasiewicz inequality is then applied,
yielding a useful estimate on the integral of the form $\int_{t_p}^{t_q} \Enorm{\dot{\vva}(t)}$.

Unfortunately, Equation~\ref{eqn:Ba13_1} is problematic. First, it is unclear
what $m(\partial V(\vva_\Gamma(t)))$ means. $V$ is a function from $\reals^N$ to $\reals$. 
So generalized gradient (or gradient for that matter) are defined in terms of elements
in $\reals^N$. In general $\vva_\Gamma$ is a shorter vector. Notwithstanding the clarification
that is needed here, the result of Equation~\ref{eqn:Ba13_1}, namely
\begin{equation}\label{eqn:Ba13_2}
\Enorm{\dot{\vva}(t)} \ge \beta\,m(\partial V(\vva(t)))
\end{equation}
cannot hold in general: Equation 8 of~\cite{Balavoine13} states that
$$
\partial V(\vva(t)) = -\vvb + \transpose{\Phi}\Phi \vva(t) + \lambda \partial C(\vva(t)),
$$
where
$\vvb = \transpose{\Phi} \vvy$ and $\vvy$ is a constant vector.
For simplicity's sake, consider the specific case when $C(a) = |a|$ (and $\beta$ corresponds to 1).
If $\vvu(t) = \vvzero$ (for example $\vvzero$ is chosen as the initial value for LCA),
$\vva(t) = \dot{\vva}(t) = \vvzero$ for at least a short period of time 
(before any of the components of $\vvu(t)$ ventures outside $[-\lambda,\lambda]$).
The left hand side of Equation~\ref{eqn:Ba13_2} is 0. On the other hand,
$$
\partial V(\vvzero) = -\vvb + \lambda \partial C(\vvzero).
$$
For any element in the set of generalized gradient, the $n$-th component is $-b_n + \gamma$
where $\gamma \in [-\lambda,\lambda]$. Thus, as long as $|b_n| = \lambda+\delta > \lambda$ 
for any one $b_n$, the norm of any element of generalized gradients is at least $\delta$.
In other words, $m(\partial V(\vvzero)) \ge \delta > 0$.
As Equation~\ref{eqn:Ba13_1} is invalid, the proof of the main result in~\cite{Balavoine13}
has a major gap.

\section*{Acknowledgment}
The author thanks Mike Davies, Tsung-han Lin and Justin Romberg
for fruitful discussions.

%\ifCLASSOPTIONcaptionsoff
%  \newpage
%\fi



%\bibliographystyle{IEEEtran}

%\bibliographystyle{IEEEtran}
\bibliographystyle{siam}
\bibliography{biblio}

\end{document}